\documentclass{article}
\usepackage[utf8]{inputenc}
\usepackage{amsmath}
\usepackage{amssymb}
\usepackage{textcomp}
\usepackage[english]{babel}
\usepackage{amsthm}
\usepackage{comment}
\usepackage{hyperref}
\usepackage{tikz}
\usetikzlibrary{arrows,decorations.pathmorphing,backgrounds,positioning,fit,matrix}
\theoremstyle{definition}
\newtheorem{definition}{Definition}[section]
\newtheorem*{remark}{Remark}
\newtheorem{theorem}{Theorem}[section]
\newtheorem{lemma}[theorem]{Lemma}
\newtheorem{proposition}[theorem]{Proposition}

\newtheorem{claim}[theorem]{Claim}
\newtheorem{fact}[theorem]{Fact}
\newtheorem*{setting}{Setting}
\makeatletter
\newcommand*{\doublerightarrow}[2]{\mathrel{
  \settowidth{\@tempdima}{$\scriptstyle#1$}
  \settowidth{\@tempdimb}{$\scriptstyle#2$}
  \ifdim\@tempdimb>\@tempdima \@tempdima=\@tempdimb\fi
  \mathop{\vcenter{
    \offinterlineskip\ialign{\hbox to\dimexpr\@tempdima+1em{##}\cr
    \rightarrowfill\cr\noalign{\kern.5ex}
    \rightarrowfill\cr}}}\limits^{\!#1}_{\!#2}}}
\newcommand*{\triplerightarrow}[1]{\mathrel{
  \settowidth{\@tempdima}{$\scriptstyle#1$}
  \mathop{\vcenter{
    \offinterlineskip\ialign{\hbox to\dimexpr\@tempdima+1em{##}\cr
    \rightarrowfill\cr\noalign{\kern.5ex}
    \rightarrowfill\cr\noalign{\kern.5ex}
    \rightarrowfill\cr}}}\limits^{\!#1}}}
\makeatother
\newtheorem*{acknowledgements}{Acknowledgements}
\title{A model theoretic proof
for o-minimal coherence theorem}
\author{Yayi Fu }
\date{}
\begin{document}
\maketitle
\begin{abstract}
    Bakker, Brunebarbe, Tsimerman showed in 
    \cite{bakker2022minimal} that the definable structure sheaf $\mathcal{O}_{\mathbb{C}^n}$ of $\mathbb{C}^n$ is a coherent $\mathcal{O}_{\mathbb{C}^n}$-module as a sheaf on the site $\underline{\mathbb{C}^n}$, 
    where the coverings are finite coverings by definable open sets. 
    In general, let $\mathcal{K}$ be an algebraically closed field of characteristic zero. 
    We give another proof of the coherence of $\mathcal{O}_{\mathcal{K}^n}$ as a sheaf of 
    $\mathcal{O}_{\mathcal{K}^n}$-modules on the site $\underline{\mathcal{K}^n}$ using spectral topology on the type space $S_n(\mathcal{K})$. 
    (Here $S_n(\mathcal{K})$ means
    $S_{2n}(\mathcal{R})$ for some real closed field $\mathcal{R}$.) 
It also gives an example of how the intuition that sheaves on the type space are the same as sheaves on the site with finite coverings (see \cite[Proposition~3.2]{edmundo2006sheaf}) can be applied.
\end{abstract}
\section{Introduction}
\indent

Let $\mathcal{O}_{\mathbb{C}^n}$ denote the sheaf of rings where $\mathcal{O}_{\mathbb{C}^n}(U)$ is the ring of holomorphic functions defined on $U$, for each $U\subseteq \mathbb{C}^n$ open. It's also an $\mathcal{O}_{\mathbb{C}^n}$-module.
\\
\indent
In complex analysis, it is well-known that
\begin{fact}\cite{oka1950fonctions}\label{oka}
    (Oka) For any positive integer $n$, $\mathcal{O}_{\mathbb{C}^n}$ is a coherent $\mathcal{O}_{\mathbb{C}^n}$-module.
    i.e.
    $\mathcal{O}_{\mathbb{C}^n}$ satisfies that
    \begin{enumerate}
        \item $\mathcal{O}_{\mathbb{C}^n}$ locally finite.
        \item Every relation sheaf of $\mathcal{O}_{\mathbb{C}^n}$is locally finite.
    \end{enumerate}
\end{fact}

This result is generalized in \cite{peterzil2008complex} to the case of any algebraically closed field $\mathcal{K}$ 
of characteristic 0.
\begin{fact}\label{ps}
     (Peterzil, Starchenko) For any positive integer $n$, $\mathcal{O}_{\mathcal{K}^n}$ is a coherent $\mathcal{O}_{\mathcal{K}^n}$-module.
\end{fact}
In fact \ref{oka}, \ref{ps},
a sheaf means the usual sheaf in e.g.
\cite[Chapter~II]{hartshorne2013algebraic}.
In \cite{bakker2022minimal}, coherence theorem is proved on the site $\underline{\mathbb{C}^n}$ where the coverings are finite coverings by definable open sets: 
\begin{fact}
     (Bakker, Brunebarbe, Tsimerman) The definable structure sheaf $\mathcal{O}_{\mathbb{C}^n}$ of $\mathbb{C}^n$
is a coherent $\mathcal{O}_{\mathbb{C}^n}$-module (as a sheaf on the site $\underline{\mathbb{C}^n}$).
\end{fact}
(The sheaves on a site in \cite{bakker2022minimal} are different from the usual sheaves defined in \cite{hartshorne2013algebraic}. 
We will explain more in later sections.)
\\\indent
In this paper, we use a method different from the one used in \cite{bakker2022minimal} to prove the coherence of $\mathcal{O}_{\mathcal{K}^n}$ as a sheaf on the site $\underline{\mathcal{K}^n}$,
     where $\mathcal{K}$ is an algebraically closed field of characteristic $0$.
\begin{theorem}\label{mainthm1}
     The definable structure sheaf $\mathcal{O}_{\mathcal{K}^n}$ of $\mathcal{K}^n$ is a coherent $\mathcal{O}_{\mathcal{K}^n}$-module as a sheaf on the site $\underline{\mathcal{K}^n}$.
\end{theorem}
Motivation of this proof comes from 
\cite[Proposition~3.2]{edmundo2006sheaf} which says we can consider a sheaf on the site $\underline{\mathcal{K}^n}$ the same as a usual sheaf in \cite[Chapter~II]{hartshorne2013algebraic} on the type space $S_n(\mathcal{K})$ with spectral topology.
\\
\indent
Section \ref{prelim} gives definitions of sites,
presheaves and sheaves on a site, spectral topology, coherence, tubular neighborhoods.
Section \ref{prf} gives the proof of theorem \ref{mainthm1}.
Section \ref{rem} shows that we can prove theorem
\ref{mainthm1} using an isomorphism of categories similar to that in \cite[Proposition~3.2]{edmundo2006sheaf}.
\begin{acknowledgements}
The author is grateful to her advisor Sergei Starchenko for the suggestion of using spectral topology and compactness to give a more model-theoretic proof
and the suggestion of using tubular neighborhoods to prove lemma \ref{tubu}.
\end{acknowledgements}
\section{Preliminaries}\label{prelim}
\subsection{Basic notions}
\begin{setting}(The same setting as in \cite{peterzil2001expansions}.)
    Let $\mathcal{K}$ be an algebraically closed field of characteristic zero. Then $\mathcal{K} = \mathcal{R}(\sqrt{-1})$ for some real closed subfield $\mathcal{R}$. 
    Such $\mathcal{R}$ is not unique.
    We fix one such $\mathcal{R}$ and fix an o-minimal expansion of the chosen real closed field. 
    The topology on $\mathcal{R}$ is generated by the definable open intervals. The topology on $\mathcal{K}$ is identified with that on $\mathcal{R}^2$. 
    When we say definable, 
    we mean definable in the o-minimal structure 
    $\mathcal{R}$ with parameters in $\mathcal{R}$.
\end{setting}
\begin{definition}
    \cite[Definition 2.1.]{peterzil2003expansions} For $U \subseteq\mathcal{K}$ 
    a definable open set and $F : U \rightarrow\mathcal{K}$ a definable function, $z_0 \in U$,
    we say that $F$ is \emph{$\mathcal{K}$-differentiable at $z_0$} if the limit as $z$ tends to $z_0$ in $\mathcal{K}$ of $(f(z) - f(z_0))/(z - z_0)$ exists in $\mathcal{K}$ (all operations taken in $\mathcal{K}$, while the limit is taken in the topology induced on $\mathcal{K}$ by $\mathcal{R}^2$).
\end{definition}
\begin{definition}
  \cite[Definition 2.8.]{peterzil2003expansions} Let $V \subseteq \mathcal{K}^n$ be a definable open set,
  $F : V \rightarrow\mathcal{K}$ a definable map.
  $F$ is called \emph{$\mathcal{K}$-differentiable} on $V$
  if it is continuous on $V$ and for every $(z_1,..., z_n) \in V$ and $i = 1,...,n$,
  the function $F(z_1,...,z_{i-1},-,z_{i+1},...,z_n)$ is $\mathcal{K}$-differentiable in the $i$-th variable at $z_i$
  (in other words, $F$ is continuous on $V$ and $\mathcal{K}$-differentiable in each variable separately).
\end{definition}
\subsection{Spectral topology}
\begin{definition}\cite[Definition~2.2.]{edmundo2006sheaf}
    Let $X \subseteq \mathcal{R}^m$ be a definable set (with parameters in $\mathcal{R}$). 
    \textit{The o-minimal spectrum $\Tilde{X}$ of $X$} is the set of 
    complete $m$-types $S_m(\mathcal{R})$ of the first 
    order theory $Th_{\mathcal{R}}(\mathcal{R})$ which imply a formula defining $X$.
    This is equipped with the topology generated by the basic open sets of 
    the form $\Tilde{U} = \{\alpha\in \Tilde{X} : U \in \alpha\}$, where $U$ is a definable, 
    relatively open subset of $X$, and $U \in\alpha$ means the formula defining $U$ is in $\alpha$.
    We call this topology on $X$ the \emph{spectral topology}. 

\end{definition} 
\subsection{Sheaves on the type space}
\indent

Let $S_n(\mathcal{K})$ denote
$S_{2n}(\mathcal{R})$.
We use this unconventional notation to emphasize that we are considering functions on $\mathcal{K}^n$.
\\
\indent
Given a definable open set $U\subseteq\mathcal{K}^n$,
let $\mathcal{O}_{\mathcal{K}^n}(\Tilde{U})$ be the ring of $\mathcal{K}$-differentiable functions defined on $U$.
It's easy to see that this defines a sheaf on the type space $S_n(\mathcal{K})$ with spectral topology. (We mean the usual notion of sheaves.)
\\
\indent
Let $\mathcal{O}_{\mathcal{K}^n}$ denote the sheaf of rings where $\mathcal{O}_{\mathcal{K}^n}(\Tilde{U})$ is the ring of $\mathcal{K}$-differentiable functions defined on $U$, for each $U\subseteq \mathcal{K}^n$ definable open.
\begin{definition}
Given $p\in S_n(\mathcal{K})$, let \emph{$\mathcal{O}_p$}
denote the set of germs for functions 
\begin{equation*}
    \text{
$\{f:U\rightarrow\mathcal{K}:U$ is some open definable set such that $p\in \Tilde{U}$}
\end{equation*}
\begin{equation*}
 \text{
and 
$f$ is $\mathcal{K}$-holomophic on $U$ $\}$. }
\end{equation*}
Given a definable set $A\subseteq\mathcal{K}^n$, 
let \emph{$\mathcal{I}_p(A)\subseteq\mathcal{O}_p$} denote the set of germs for functions 
\begin{equation*}
\text{
$\{f:U\rightarrow\mathcal{K}:U$ is some open definable set such that $p\in \Tilde{U}$, }
\end{equation*}
\begin{equation*}
    \text{
$f$ is $\mathcal{K}$-holomophic on $U$ and $\forall x\in A\cap U$, $f(x)=0$ $\}$. }
\end{equation*}
Let $g_1,...,g_t\in\mathcal{O}_p$. 
Let \emph{$R_p(g_1,...,g_t)$} denote the set
\begin{equation*}
   \text{$\{(f_1,...,f_t)\in\mathcal{O}^t_p:f_1g_1+...+f_tg_t=0\}$.   
 }
\end{equation*}
\end{definition}

\indent

Notice that for a sheaf on the type space $S_n(\mathcal{K})$, we mean a usual sheaf as
defined in \cite[Chapter~2]{hartshorne2013algebraic}.
In the next subsection, we define a different notion of sheaves, the sheaves on a site,
where coverings are finite.
\subsection{o-minimal site}
\indent

We translate definitions about sites in
\cite{stacks-project} into o-minimal context: (For the formal definitions related to sites and how the usual notion of sheaves defined in \cite[Chapter~2]{hartshorne2013algebraic} are defined in the context of sites,
see \cite[Part~1, Chapter~7]{stacks-project}.)
\begin{definition}\label{site}
      \cite[Part~1, Chapter~7, Definition~6.2]{stacks-project} Let $X\subseteq\mathcal{K}^n$ be a definable set. The \emph{o-minimal site} $\underline{X}$ on $X$ consists of definable (relative) open subsets of $X$, together with $Cov(X):=\{(U,\{U_i\}_{i=1}^k):U,U_1,...,U_k\subseteq X$ definable open, $\{U_i\}_{i=1}^k$ a finite covering of $U$ $\}$. 
      (This means, in the formal definition of a site,
      the objects of the category are definable open subsets of $X$;
      the morphisms of the category are inclusions;
      the coverings are finite coverings by definable open sets.)
 \end{definition}
 \begin{definition}\label{presheaf}
     \cite[Part~1,
     Chapter~6,
     Section~5]{stacks-project} A \emph{presheaf} of abelian groups 
     (resp. rings) on an o-minimal site $\underline{X}$ is defined the same as usual:\\
\indent     Let $X$ be a topological space. A \emph{presheaf} $\mathcal{F}$ of abelian groups (resp. rings)
     on an o-minimal site $\underline{X}$ consists of the
following data:
\begin{enumerate}
    \item [(a)]
 a collection of non empty abelian groups (resp. rings) $\mathcal{F}(U)$ associated with every definable open set $U \subseteq X$,
 \item[(b)] 
a collection of morphisms of abelian groups (resp. rings) $\rho_{U,V} : \mathcal{F}(V)\rightarrow \mathcal{F}(U)$ defined whenever $U \subseteq V$ and satisfying
the transitivity property,
\item[(c)] 
 $\rho_{U,V} \circ \rho_{V,W} = \rho_{U,W}$ for $U \subseteq V \subseteq W$, $\rho_{U,U} = Id_U$ for every $U$. 
\end{enumerate}
 \end{definition}
 \begin{definition}\label{presheaf}
     \cite[Part~1,
     Chapter~6,
     Definition~6.1]{stacks-project} 
     Let $X$ be a topological space.
     Let 
     $\mathcal{O}$ be a presheaf of rings on the o-minimal site $\underline{X}$.
     A \textit{presheaf of $\mathcal{O}$-modules $\mathcal{F}$} 
     on an o-minimal site $\underline{X}$ 
     is a presheaf $\mathcal{F}$ of abelian groups with the
following additional data:
\begin{enumerate}
    \item [(a)]   For every definable open set $U \subseteq X$,
    $\mathcal{F}(U)$ is a
    non empty $\mathcal{O}(U)$-module;
 \item[(b)] 
for every definable open $U \subseteq X$ the $\mathcal{O}(U)$-module structure of
$\mathcal{F}(U)$ is
compatible with restriction mappings (of $\mathcal{F}$ and $\mathcal{O}$). 
i.e. for definable open $U\subseteq V\subseteq X$, 
$r\in\mathcal{O}(V)$, $x\in \mathcal{F}(V)$,
$\rho_{U,V}(r)\tau_{U,V}(x)= \tau_{U,V}(rx)$.
 
\end{enumerate}
 \end{definition}
 \begin{definition}\label{sheafonsite}
       \cite[Part~1, Chapter~6, Definition~7.1.]{stacks-project} Let $\underline{X}$ be an o-minimal site, and let $\mathcal{F}$ be a presheaf of abelian groups (resp. rings, $\mathcal{O}$-modules) on $\underline{X}$. We say $\mathcal{F}$ is a \emph{sheaf} if for every definable open $U\subseteq X$ and every definable open finite covering $\{U_i\}_{i=1}^k$ of $U$, 
       \begin{enumerate}
           \item [(i)]
        if $(s_i)_{i=1}^k$ satisfies $s_i\in\mathcal{F}(U_i)$ for each $i$ and $s_i|_{U_i\cap U_j}=s_j|_{U_i\cap U_j}$ for each pair $i,j$, then there is a unique $s\in U$ such that $s|_{U_i}=s_i$ for each $i$; 
        \item[(ii)] 
        for $s,t\in\mathcal{F}(U)$, if $s|_{U_i}=t|_{U_i}$ for each $i$ then $s=t$.
       \end{enumerate}
 \end{definition}
 \begin{definition}\label{injsur}
    \cite[Part~1,
    Chapter~7.
    Definition~11.1.]{stacks-project} Let  $\underline{X}$ be an o-minimal site, and let $\varphi: \mathcal{F}\rightarrow \mathcal{G}$ be a map of sheaves of modules. (i.e. $\varphi$ is a morphism of $\mathcal{F}$ and $\mathcal{G}$ considered as presheaves.
    A presheaf morphism is, as usual, a map compatible with the restiction maps.)
    \begin{enumerate}
        \item [(1)]
    We say that $\varphi$ is \emph{injective} if for every definable open $U\subseteq X$ the map $\varphi: \mathcal{F}(U) \rightarrow\mathcal{G}(U)$ is injective.
\item[(2)] 
We say that $\varphi$ is \emph{surjective} if for every definable open $U\subseteq X$ and every section $s \in \mathcal{G}(U)$ there exists a finite covering $\{U_i\}_{i=1}^k$ of $U$ such that for each $i$, $U_i$ is definable open and the restriction
$s|_{U_i}$ is in the image of $\varphi: \mathcal{F}(U_i) \rightarrow \mathcal{G}(U_i)$.
\end{enumerate}
 \end{definition}
 \begin{definition}\label{coherencedef} (\cite[Definition 2.13]{bakker2022minimal})  
 Let $\underline{X}$ be an o-minimal site.
 Given an $\mathcal{O}_X$-module $M$,
 we say that $M$ is of \emph{finite type} (as an $\mathcal{O}_X$-module) if there exists
a finite definable open (relative to $X$) cover $X_i$ of $X$ and surjections $\mathcal{O}^n_{X_i} \twoheadrightarrow M_{X_i}$ for some positive
integer $n$ on each of those open sets. 
We say that $M$ is \emph{coherent} (as an $\mathcal{O}_X$-module) if it is of finite type, 
and given any
definable open $U \subseteq X$ and any $\mathcal{O}_U$ -module homomorphism $\varphi : \mathcal{O}^n_U \rightarrow M_U$ ,
the kernel of $\varphi$ is of finite type.
 \end{definition}
 \begin{remark}
 By definition \ref{injsur} and definition \ref{coherencedef}, given a definable open $U$ and an $\mathcal{O}_U$-module $M$, to show that $M$ is of finite type, it suffices to show that there exist a finite family of definable open sets $U_1,...,U_k$ covering $U$ and sheaf morphisms $\varphi_i:\mathcal{O}_{U_i}\rightarrow M_{U_i}$, $i=1,...,k$ such that for any fixed $i$, for any definable open $V\subseteq U_i$ and every section $s\in M(V)$, there exist a finite family of definable open sets $V_1,...,V_l$ covering $V$ and for each $j\in\{1,...,l\}$, $t_{j}\in \mathcal{O}_{V_j}$ such that $\varphi(V_j)(t_j)=s|_{V_j}$.  
 \end{remark}
 \subsection{Motivation}
 Let $X\subseteq\mathcal{K}^n$ be a definable set.
 \begin{definition}\cite[Definition~2.2]{edmundo2006sheaf}
     For the o-minimal spectrum $\Tilde{X}$ of $X$, since it is a 
    topological space, we use the classical notation $Sh(\Tilde{X})$ 
    to denote the category of sheaves of abelian groups on $\Tilde{X}$.
Since the topology on the o-minimal spectrum $\Tilde{X}$ of $X$ is generated by the constructible open subsets, i.e. sets of the form $\Tilde{U}$ with $U$ an open definable subset of $X$, a sheaf on $\Tilde{X}$ is determined by its values on the sets $\Tilde{U}$ with $U$ an open definable subset of $X$. (We may also consider $\Tilde{X}$ as a site where
the objects of the category are the $\Tilde{U}$'s where each $U$ is some definable open subset of $X$;
      the morphisms of the category are inclusions;
      the coverings are \textit{any} coverings by the $\Tilde{U}$'s.)
 \end{definition}
\begin{definition}\cite[Definition~3.1.]{edmundo2006sheaf}
    We denote by $Sh_{dtop}(X)$ the category of sheaves of abelian groups on $X$ with respect to the o-minimal site on $X$. 
    \\
    Thus, for a definable set $X$, we define the functor of the categories of sheaves of abelian groups
$Sh_{dtop}(X) \rightarrow Sh(\Tilde{X})$
which sends $F \in Sh_{dtop}(X)$ into $\Tilde{F}$ where, 
for $U$ an open definable subset of $X$, 
we
define $\Tilde{F}(\Tilde{U}) = \{\Tilde{s}: s \in F(U)\} \simeq F(U)$,
and
$Sh(\Tilde{X}) \rightarrow Sh_{dtop}(X)$
which sends $\Tilde{F}$ into $F$ where, 
for $U$ an open definable subset of $X$, 
we define $F(U) = \{ s : \Tilde{s} \in \Tilde{F}(\Tilde{U})\} \simeq \Tilde{F} (\Tilde{U})$.
\end{definition}
The following fact is the motivation for our proof in section \ref{prf}.
It says that a sheaf on the site $\underline{\mathcal{K}^n}$
is the same as a usual sheaf in \cite[Chapter~II]{hartshorne2013algebraic} on the type space $S_n(\mathcal{K})$ with spectral topology.
\begin{fact}
    \cite[Proposition~3.2]{edmundo2006sheaf}
    $Sh(\Tilde{X})$ and $Sh_{dtop}(X)$ are isomorphic.
\end{fact}
\subsection{Tubular neighborhood}
\indent

\cite[Theorem~2.56]{peterzil2001expansions}
roughly says that given a $\mathcal{K}$-holomorphic function $f$ and $p\in\mathcal{K}^n$,
the number of zeroes is fixed locally around $p$.
In this section,
we prove the following lemma, 
which says that for all $p\in S_n(\mathcal{K})$,
the number of zeroes is fixed locally.
This lemma will be used in the proof of the type version of Weierstrass division theorem.
\\
\indent
Let $\pi:\mathcal{K}^n\rightarrow \mathcal{K}^{n-1}$,
$\pi_n:\mathcal{K}^n\rightarrow \mathcal{K}$
denote the projection onto the first $(n-1)$
coordinates and the projection onto the $n$-th
coordinate resp.
\\
\indent
For notational convenience, 
given $p\in S_n(\mathcal{K})$ and a definable set $U$,
when we say \textquotedblleft$p\in U$",
we actually mean \textquotedblleft$p\in\Tilde{U}$".
Similarly, when we say \textquotedblleft an open neighborhood $U$
of $p$",
we actually mean \textquotedblleft an open neighborhood $\Tilde{U}$
of $p$".
\begin{lemma}\label{tubu}
  Let $p\in S_n(\mathcal{K})$.
Fix $f\in \mathcal{O}_p$ and an open neighborhood $U$ of $p$ on which $f$
is defined and is $\mathcal{K}$-differentiable.
 Suppose for all $y\in\pi(U)$,
 there are finitely many zeroes of 
 $f(y,-)$ in $U_y:=
 \{x\in\mathcal{K}: (y,x)\in U\}$, 
 counting multiplicity.
 \\
 \indent
 Then there exist $i\in\mathbb{N}$ and
 $V\subseteq U$ a definable open neighborhood of $p$ 
        such that for any $y\in \pi(V)$,
        there are exactly $i$ zeroes of
        $f(y,-)$ in $V_y$ counting multiplicity.
\end{lemma}
We need some basic definitions and facts about 
o-minimal structures.
\begin{definition}\cite[Chapter~3]{van1998tame}
    Call a set $Y\subseteq \mathcal{R}^{m+1}$ is \emph{finite over $\mathcal{R}^m$} if for each $x\in \mathcal{R}^m$ the fiber $Y_x:= \{r\in\mathcal{R}:
    (x,r)\in Y\}$ is finite;
    call $Y$ \emph{uniformly finite over $\mathcal{R}^m$}
    if there is  $N\in \mathbb{N}$ such that $|Y_x| \leq N $ for all $x \in \mathcal{R}^m$.
\end{definition}
\begin{fact}\label{finiteness}
    \cite[Chapter~3, Lemma (2.13)]{van1998tame}
 (UNIFORM FINITENESS PROPERTY). Suppose the definable subset $Y$ of $\mathcal{R}^{m+1}$ is finite over $\mathcal{R}^{m}$. 
 Then $Y$ is uniformly finite over $\mathcal{R}^{m}$.
\end{fact}
\begin{fact}\label{numofzeroes}
    \cite[Theorem~2.56.]{peterzil2001expansions}
    Let $W \subseteq \mathcal{R}^n$, 
    $U \subseteq \mathcal{K}$
    be definable open sets,
    $F : U \times W \rightarrow \mathcal{K}$ a definable continuous function such that for every 
    $w  \in W$, $F(-,w )$ is a $\mathcal{K}$- differentiable function on $U$.
Take $(z_0,w_0) \in U \times W$ and suppose that $z_0$ is a zero of order $m$ of
$F(-,w_0)$.
\\
\indent
Then for every definable neighborhood $V$ of $z_0$ there are definable open neighborhoods $U_1 \subseteq V$ of $z_0$ and
$W_1 \subseteq W$ of $w_0$ such that
$F(-,w)$ has exactly $m$ zeroes in $U_1$ (counted with multiplicity) for every 
$w \in W_1$.
\end{fact}
We may assume in fact \ref{numofzeroes} that
$U_1, W_1$ are open balls:
Let $ U_1\subseteq U, W_1\subseteq W$ 
be definable neighborhoods of $z_0$,
$w_0$ resp. such that for all $w\in W_1$,
$F(-,w)$ has exactly $m$ zeroes in $U_1$.
Let $U_2\subseteq U_1$ and $W_2\subseteq W_1$ be definable open balls of $z_0$, $w_0$ resp. 
Then for all $w\in W_2$,
$F(-,w)$ has $\leq m$ zeroes in $U_2$.
By fact \ref{numofzeroes},
there exist $U_3\subseteq U_2$,
$W_3\subseteq W_2$ definable open neighborhoods of $z_0$, $w_0$ resp.
such that 
$F(-,w)$ has exactly $m$ zeroes in $U_3$ (counted with multiplicity) for every 
$w \in W_3$. 
Let $W_4\subseteq W_3$ be a definable open ball around $w_0$.
Then for all $w\in W_4$,
$F(-,w)$ has $\geq m$ zeroes in $U_2$.
Hence $U_2$ and $W_4$ are open balls satisfying the conclusion of fact \ref{numofzeroes}.
\begin{definition}
    \cite[Section~6.2]{coste2000introduction}
    A \emph{$C^k$ cylindrical definable cell decomposition of $\mathcal{R}^n$} is a cdcd satisfying extra smoothness conditions which imply,
    in particular, that each cell is a $C^k$ submanifold of $\mathcal{R}^n$.
    \begin{itemize}
        \item 
 A $C^k$ cdcd of $\mathcal{R}$ is any cdcd of $\mathcal{R}$ 
 (i.e. a finite subdivision of $\mathcal{R}$).
\item  If $n>1$, 
a $C^k$ cdcd of $\mathcal{R}^n$ is given by a $C^k$ cdcd of $\mathcal{R}^{n-1}$ and,
for each cell
$D$ of $\mathcal{R}^{n-1}$,
definable functions of class $C^k$
$\zeta_{D,1} <...<\zeta_{D,l(D)}
:D\rightarrow \mathcal{R}$.
The cells of $\mathcal{R}^n$ are, of course, the graphs of the $\zeta_{D,i}$ and the bands delimited by these graphs.
 \end{itemize}
\end{definition}
\begin{fact}\label{celldecomp}
    \cite[Theorem~6.6]{coste2000introduction} 
    ($C^k$ Cell Decomposition: $C^k$ $CDCD_n$)
    Given finitely many definable subsets $X_1,...,X_l$ of $\mathcal{R}^n$,
    there is a $C^k$ cdcd of $\mathcal{R}^n$ adapted to $X_1,...,X_l$
    (i.e. each $X_i$ is a union of cells).
\end{fact}
\begin{fact}\label{pwcts}
\cite[Theorem~6.7]{coste2000introduction}
(Piecewise $C^k$ $PC^k_n$)
Given a definable function $f:A\rightarrow\mathcal{R}^n$,
where $A$ is a definable subset of $\mathcal{R}^n$,
there is a finite partition of $A$ into definable $C^k$ submanifolds 
$\mathcal{C}_1,...,\mathcal{C}_l$, such that each restriction $f|\mathcal{C}_i$ is $C^k$.
\end{fact}
\begin{definition}
\cite[Chapter~6]{coste2000introduction}
Let $M \subseteq \mathcal{R}^n$ be a definable $C^k$ submanifold 
(we always assume $1\leq  k <\infty$).
The \emph{tangent bundle $TM$} is the set of $(x,v) \in M \times \mathcal{R}^n$ such that $v$ is a tangent vector to $M$ at $x$.
    The \emph{normal bundle $NM$} is the set of $(x,v)$ in $M \times\mathcal{R}^n$ such that $v$ is orthogonal to $T_xM$.
    This is a $C^{k-1}$ submanifold of
    $\mathcal{R}^n \times\mathcal{R}^n$,
    and it is definable since $TM$ is definable.
\end{definition}
Let $\varphi$ be the function 
$\varphi:N\mathcal{D}\rightarrow \mathcal{R}^{2n-2}$
where
$\varphi(x,v)=x+v$. 
\begin{fact}
    \cite[Theorem~6.11]{coste2000introduction}
    (Definable Tubular Neighborhood) Let $M$ be a definable $C^k$ submanifold of $\mathcal{R}^n$.
    There exists a definable open neighborhood $U$ of the zero-section
    $M\times \{0\}$ in the normal bundle $NM$ such that the restriction $\varphi|U$ is a $C^{k-1}$ diffeomorphism onto an open neighborhood $\Omega$ of $M$ in $\mathcal{R}^n$.
    Moreover, we can take $U$ of the form
\begin{equation*}
    U =\{(x,v)\in NM : \|v\|<\epsilon(x)\}, 
\end{equation*}
where $\epsilon$ is a positive definable $C^k$ function on $M$.
\end{fact}
\begin{fact}\label{lowerbd}
   \cite[Lemma~6.12]{coste2000introduction} 
   Let $M$ be a definable $C^k$ submanifold of $\mathcal{R}^n$,
   closed in $\mathcal{R}^n$.
   Let $\psi : M \rightarrow \mathcal{R}$ be a positive definable function, 
   which is locally bounded from below by positive constants
   (for every $x$ in $M$,
   there exist $c > 0$ and a neighborhood $V$ of $x$ in $M$ such that $\psi>c$ on $V$). 
   Then there exists a positive definable
   $C^k$ function $\epsilon:M \rightarrow
   \mathcal{R}$ such that $\epsilon<\psi$ on $M$.
\end{fact}
\begin{fact}\label{tubularsp}
    \cite[Lemma~6.15] {coste2000introduction}
    The definable tubular neighborhood
    in \cite[Theorem~6.11]{coste2000introduction} holds if $M$ is closed in $\mathcal{R}^n$, or definably $C^k$ diffeomorphic to a closed submanifold in some $\mathcal{R}^m$ (e.g. if $M$ is a cell of a $C^k$ cdcd).
\end{fact}
Now we prove lemma \ref{tubu} by imitating the proof of fact \ref{tubularsp}.
\begin{proof}
Since the set $Z(f):=\{z: f(z)=0\}$
is a closed set, 
if $p\notin Z(f)$,
then $U\setminus Z(f)$ is a
definable open set satisfying the lemma.
Hence,
may assume $p\in Z(f)$.
\\
\indent
Let $p\in S_n(\mathcal{K})$.
Fix $f\in \mathcal{O}_p$ and an open neighborhood $U$ of $p$ on which $f$
is defined and is $\mathcal{K}$-differentiable.
\\
\indent
By fact \ref{finiteness}, 
since for all $y\in\pi(U)$,
 there are finitely many zeroes of 
 $f(y,-)$ in $U_y$, 
there is $m$ such that for all $y\in\pi(U)$,
 there are $\leq m$ many zeroes of 
 $f(y,-)$ in $U_y$.
\\
\indent
Let $M_i$ be the set 
\begin{equation*}
\{x\in U: \pi_n(x) \text{ is a zero of order $i$  for the function } f(\pi(x), -)\}.
\end{equation*}
Then $p\in M_i$ for some $i$.
Fix such $i$ and
let $M=M_i$.
Let $\pi(M)$ denote $\{y\in \mathcal{K}^{n-1}:
\exists z$ $(y,z)\in M\}$.
\\
\indent
Define $F_1,...,F_m:\pi(M)\rightarrow \mathcal{K}$ as follows:
\begin{equation*}
    F_1(u)=\text{ the least } v 
    \text{ such that }(u,v)\in M.
\end{equation*}
Suppose $F_1,...,F_j$ are defined. 
Let
\begin{equation*}
    F_{j+1}(u)=\text{ the least } v
    \text{ such that }
(u,v)\in M \text{ and }
v\notin\{F_1(u),...,F_j(u)\}
\end{equation*}
if such $v$ exists;
otherwise, let
\begin{equation*}
    F_{j+1}(u)=F_j(u).
\end{equation*}
Then $M=(\pi(M)\times\mathcal{K})\cap M=
\underset{j=1}{\overset{m}{\bigcup}}
\{(u,F_j(u)):u\in\pi(M)\}$.
\\
Take $j$ such that $p\in\{(u,F_j(u)):
u\in\pi(M)\}$.
\\
\indent
Let $F: \pi(M)\rightarrow M$ be the definable function $F_j$.
Then for each $y\in\pi(M)$,
$F(y)$ is a zero of order $i$ for the function $f(y,-)$ and $p\in grpah(F)$.
Then there is a $C^k$-cell $\mathcal{C}\subseteq \pi(M)$ such that 
$F$ is continuous on $\mathcal{C}$ and
$p\in (\mathcal{C}\times\mathcal{K})\cap
M$.
(\cite{coste2000introduction}
the theorem for cell decomposition and the theorem for piecewise continuous)
\\
\indent
Recall that $\sup $ and $\inf$ 
exist for definable sets by
\cite[Chapter~1, Lemma~(3.3)~(i)]{van1998tame}.
Also recall fact \ref{numofzeroes},
which roughly says that the number of zeroes remains the same in a small neighborhood around a fixed zero.
Then we can define functions $\alpha$, $\beta$
as follows:
\\
\indent
Define $\alpha:\mathcal{C}\rightarrow\mathcal{R}^{>0}$ by 
\begin{equation*}
    \alpha(x)=
    \frac{1}{2}
   \sup\{r\in \mathcal{R}:\text{ there is $s>0$ such that for all $y\in B(x,r)$}
   \end{equation*}
   \begin{equation*}
      \text{ there are exactly $i$ zeroes in $B(F(x),s)$ for $f(y,-)\}$}.
\end{equation*}
Define $\beta:\mathcal{C}\rightarrow\mathcal{R}^{>0}$ by
\begin{equation*}
    \beta(x)=
   \sup\{r\in \mathcal{R}:\text{ for all $y\in B(x,\alpha(x))$ there are exactly $i$ zeroes }
   \end{equation*}
   \begin{equation*}
       \text{ in $B(F(x),r)$ for $f(y,-)\}$.}
\end{equation*}
Then there is a $C^k$-cell
$\mathcal{D}\subseteq\mathcal{C}$
such that $\alpha, \beta$ are continuous on $\mathcal{D}$.
Note: $\beta(x)$ satisfies that
for all $y\in B(x,\alpha(x))$ there are exactly $i$ zeroes in $B(F(x),\beta(x))$ for $f(y,-)\}$ 
because
$B(F(x),\beta(x))=\underset{r<\beta(x)}{\bigcup} B(F(x),r)$.
\\
\indent
As in \cite{coste2000introduction},
let $\varphi$ be the function 
$\varphi:N\mathcal{D}\rightarrow \mathcal{R}^{2n-2}$
where
$\varphi(x,v)=x+v$ and
let $Z$ be the subset of $(x, v)$ in $N\mathcal{D}$ such that
\begin{equation*}
d_{(x,v)}\varphi : 
T_{(x,v)}(N\mathcal{D}) \rightarrow \mathcal{R}^{2n-2}
\end{equation*}
is not an isomorphism. 
Define $\theta:\mathcal{D}\rightarrow\mathcal{R}^{>0}$ such that
\begin{enumerate}
    \item $\theta(x)\leq \min\{1, dist ((x,0),Z),
    \alpha(x)\}$
    \item 
    $\theta(x)\leq 
    \inf\{r\in \mathcal{R}:
        \exists(y,w)\in N\mathcal{D}$ $
        \exists v\in N_x\mathcal{D}$ 
       $\| w\|\leq 
       \| v\|=r$ and 
       $y+w=x+v$
        $\}$.
\end{enumerate}
By continuity of $\alpha$ and by the proof of fact \ref{tubularsp},
$\theta$ is locally bounded below.
 By fact \ref{lowerbd},
 there is a definable continuous $\epsilon:\mathcal{D}
    \rightarrow \mathcal{R}^{>0}$
    such that $\epsilon<\theta$.
    (We are not using fact 
    \ref{pwcts} to get a cell $\mathcal{D}'\subseteq \mathcal{D}$
    on which $\theta$ is continuous since
    $N\mathcal{D}'$ might be very 
    different from $N\mathcal{D}$.)
    Define a set 
    \begin{equation*}
        V=
        \{(y,z)\in\mathcal{K}^{n-1}\times
        \mathcal{K}:
        y=u+v \text{ for some unique }
        u\in\mathcal{D}, v\in N_u\mathcal{D}
        \text{ and } 
        \end{equation*}
        \begin{equation*}
                  \|v\|<\epsilon(u),
        z\in B(F(u),\beta(u))
        \}
    \end{equation*}
    By the proof of fact \ref{tubularsp},
     $\pi(V)$ is definable open and
        $\varphi|\pi(V)$ is a diffeomorphism.
    \begin{claim}
        $V$
        is a open neighborhood of $p$ 
        such that for any $y\in \pi(V)$,
        there are exactly $i$ zeroes of
        $f(y,-)$ in $V_y$.
    \end{claim}
    \begin{proof}
        $p\in V$ since 
        $p\in graph(F)\subseteq
        V$.
        \\
        \indent
        Fix $y\in\pi(V)$. 
        Let $u\in\mathcal{D}$,
        $v\in N_u\mathcal{D}$ be the unique 
        elements such that $y=u+v$.
        Then $z\in V_y$ iff
        $z\in B(F(u),\beta(u)))$.
        By the choice of $\beta$, $\epsilon$,
        there are exactly $i$ zeroes of 
        $f(y,-)$
        in $B(F(u),\beta(u))$.
        \\
        \indent
        We now show 
        that
        $V$ is open.
    Let $(y,z)\in V$. 
    Write $y=u+v$ where $u\in\mathcal{D}$ and
        $v\in N_u\mathcal{D}$ are unique.
        By continuity of $\alpha,\beta,F$,
        there is $s>0$ such that
        for all $u'\in\mathcal{D}$ with $\|u'-u\|<s$,
        \begin{itemize}
            \item 
        $\alpha(u')>
       \frac{1}{2}(\alpha(u)-\|v\|)+\|v\|$,
       \item 
        $\epsilon(u')>
       \frac{1}{2}(\epsilon(u)-\|v\|)+\|v\|$,
       \item 
        $\beta(u')>\|z-F(u)\|+\frac{1}{2}(\beta(u)-\|z-F(u)\|)$
        and
        \item 
        $\|F(u)-F(u')\|\leq
        \frac{1}{4}(\beta(u)-\|z-F(u)\|)$.
         \end{itemize}
        Since by the proof of fact \ref{tubularsp},
        $\varphi|\pi(V)$ is a diffeomorphism,
        there is $r>0$ such that for
        all $y'\in \mathcal{K}^{n-1}$
        with $\|y'-y\|<r$,
        $y'=u'+v'$ where $u'\in\mathcal{D}$,
        $v'\in N_{u'}\mathcal{D}$ are unique
    and $\|u'-u\|<s$, 
    $\|v'-v\|<
    \min\{\frac{1}{4}(\epsilon(u)-\|v\|),
    \frac{1}{4}(\alpha(u)-\|v\|)\}$.
    Then for $(y',z')\in
    B(y,r)\times B(z, \frac{1}{4}(\beta(u)-\|z-F(u)\|))$,
    write $y'$ as $y'=u'+v'$ with
    unique $u'\in\mathcal{D}$ and
        $v'\in N_u\mathcal{D}$.
        We have
    \begin{equation*}
        \|v'\|\leq
    \|v'-v\|+\|v\|\leq 
  \min\{\frac{1}{4}(\alpha(u)-\|v\|)+\|v\|,
  \frac{1}{4}(\epsilon(u)-\|v\|)+\|v\|\}
  \end{equation*}
  \begin{equation*}
    <\min\{\epsilon(u'),\alpha(u')\}
    \end{equation*}
    and
    \begin{equation*}
    \|z'-F(u')\|\leq
    \|z'-z\|+\|z-F(u)\|+\|F(u)-F(u')\|\leq
    \end{equation*}
    \begin{equation*}
        \frac{1}{4}(\beta(u)-\|z-F(u)\|)+
        \|z-F(u)\|+
        \frac{1}{4}(\beta(u)-\|z-F(u)\|)\leq
    \end{equation*}
    \begin{equation*}
        \|z-F(u)\|+
        \frac{1}{2}(\beta(u)-\|z-F(u)\|)
        <\beta(u')
    \end{equation*}
    Hence $B(y,r)\times B(z, \frac{1}{4}(\beta(u)-\|z-F(u)\|))\subseteq V$
    and $V$ is open.
    \end{proof}
      
\end{proof}
\section{Proof}\label{prf}
\textit{Outline of the proof:}
Given $p\in S_n(\mathcal{K})$,
use lemma \ref{tubu} to get a neighborhood
of $p$ with fixed number of zeroes.
Then follow the proof of \cite[Theorem~2.23]{peterzil2003expansions} to get the type version of Weierstrass division theorem.
Then the rest is just the same as in \cite{peterzil2008complex}.
\begin{theorem}
    (type version of \cite[Theorem~2.23.]{peterzil2003expansions}) \label{tpvwdiv}
 \\
Let $p\in S_n(\mathcal{K})$,
$U$ a definable open neighborhood of $p$.
Let $f(z_1,...,z_{n-1},y), g(z,y)\in\mathcal{O}_{p,n}$ be defined and $\mathcal{K}$-differentiable on $U$. 
Suppose for all $y\in\pi(U)$,
 there are finitely many zeroes of 
 $f(y,-)$ in $U_y:=
 \{x\in\mathcal{K}: (y,x)\in U\}$, 
 counting multiplicity.
\\
\indent
Then there is $k\in\mathbb{N}$, a definable open set $V\subseteq U$ with $p\in V$ and 
unique $q(z,y)\in\mathcal{O}_{V,n}$, 
$R_0(z),...,R_{k-1}(z)\in\mathcal{O}_{V,n-1}$ 
such that 
\begin{equation*}
   g(z,y)=q(z,y)f(z,y)+R_{k-1}(z)y^{k-1}+...+R_1(z)y+R_0(z)
   \text{ on $V$}.
\end{equation*}
\end{theorem}
 \begin{proof}
      Fix $g\in\mathcal{O}_{p,n}$. 
      Suppose $f,g$ are defined
      and $\mathcal{K}$- differentiable on a definable open set $U'\subseteq U$.
      By lemma \ref{tubu}, 
      there exist $k\in\mathbb{N}$ and
 $V\subseteq U'$ a definable open neighborhood of $p$ 
        such that for any $y\in \pi(V)$,
        there are exactly $k$ zeroes of
        $f(y,-)$ in $V_y$ (counting multiplicity).
         Then the rest of the proof is the same as in \cite[Theorem~2.23.]{peterzil2003expansions}.\end{proof}
 \begin{theorem}\label{tpvimmstp}
     (type version of \cite[Theorem~11.2.]{peterzil2008complex}) Assume that $U\subseteq\mathcal{K}^n$ is a definable open set and $A\subseteq U$ an irreducible $\mathcal{K}$-analytic subset of $U$ of dimension $d$. Assume also:
     \begin{enumerate}
         \item [(i)]
      The projection $\pi$ of $A$ on the first $d$ coordinates is definably proper over its image, and $\pi(A)$ is open in $\mathcal{K}^d$.
      \item[(ii)] 
      There is a definable set $S\subseteq\mathcal{K}^d$, of $\mathcal{R}$-dimension $\leq2d-2$ and a natural number $m$, such that $\pi|A$ is $m$-to-$1$ outside the set $A\cap\pi^{-1}(S)$, $\pi$ is a local homeomorphism outside of the set $\pi^{-1}(S)$, and $A\setminus\pi^{-1}(S)$ is dense in $A$.
     \item[(iii)] 
      The coordinate function $z\mapsto z_{d+1}$ is injective on $A\cap\pi^{-1}(x')$ for every $\mathcal{R}$-generic $x'\in\pi(A)$. Namely, for all $z, w\in\pi^{-1}(x')$, if $z_{d+1}=w_{d+1}$  then $z=w$.
     \end{enumerate}
     \indent
     
     Then, there is a definable open set $U'\subseteq U$ containing $A$, a natural number $s$ and $\mathcal{K}$-holomorphic functions $G_1,...,G_r,D:U'\rightarrow\mathcal{K}$, such that for every $p\in S_n(\mathcal{K})$ with $p\in A$ and $f\in\mathcal{O}_p$, if $g_1,...,g_r,\delta$ are the germs at $p$ of $G_1,...,G_r,D$, resp, then:\\
     \begin{equation} 
     f\in\mathcal{I}(A)_p\iff\exists f_1,...,f_r\in\mathcal{O}_p(\delta^sf=f_1g_1+...+f_rg_r)
     \end{equation}     
     \end{theorem}
     \begin{proof}
      Define $P_{d+1},...,P_n$ and $\{D(z')z_i-R_i(z',z_{d+1}):i=d+2,...,n\}$ satisfying \cite[Claim~11.4., Claim~11.5., Claim~11.6.]{peterzil2008complex} as in the proof of \cite[Theorem~11.2.]{peterzil2008complex}. These are $\mathcal{K}$-holomorphic functions defined on the open set $\pi(A)\times\mathcal{K}^{n-d}\supseteq A$. Hence $p\in \pi(A)\times\mathcal{K}^{n-d}$ and we can consider the germs $p_{d+1},...,p_{n}$ for $P_{d+1},...,P_{n}$ respectively, the germ $\delta$ for $D(z')$, and $r_i(z',z_{d+1})$ for $R_i(z',z_{d+1})$ in the ring $\mathcal{O}_p$. \\
      \indent
      Let $J_p$ be the ideal of $\mathcal{O}_p$ generated by the germs of $P_{d+1},...,P_n$ and $D(z')z_i-R_i(z',z_{d+1})$, $i=d+2,...,n$ at $p$. Let $s=(m-1)(n-(d+1))$. As in \cite[Theorem~11.2.]{peterzil2008complex}, for all $f\in\mathcal{O}_p$, $\delta^s f\in J_p$.\\
      \indent
      Now we show that for all $f\in\mathcal{O}_p$, if $\delta^s f\in J_p$ then $f\in\mathcal{I}(A)_p$. We follow the proof in \cite[Theorem~11.2.]{peterzil2008complex}
      \begin{claim}(type version of \cite[Claim~11.7.]{peterzil2008complex})
      \label{claimdivisible}
          For every $h(z',z_{d+1})\in\mathcal{O}_p$, if $h\in\mathcal{I}(A)_p$ then $h$ is divisible by $p_{d+1}(z',z_{d+1})$.
      \end{claim}
      \begin{proof}
      Let $h$ be defined on some definable open set $W$ with $p\in W$. 
      Since $A\subseteq\pi(A)\times\mathcal{K}^{n-d}$ and $p\in A$, may assume 
      $\pi(W)\subseteq\pi(A)$. $\forall z\in W$, if $p_{d+1}(z)=0$ then $h(z)=0$: 
      Suppose $z\in W$ and $p_{d+1}(z)=0$.
      Write $z$ as $z=(\pi(z),z_{d+1},...,z_n)$.
        If $\pi(z)\notin S$, 
       then 
      $z_{d+1}=\phi_{i,d+1}(\pi(z))$ for some $i\in\{1,...,m\}$. (The $\phi_i=
      (\phi_{i,d+1},...,\phi_{i,n})$'s are defined as in \cite[Claim~11.4.]{peterzil2008complex}.) 
      Since  
      $(\pi(z),\phi_i(\pi(z)))\in A$ 
      and $h\in \mathcal{I}(A)_p$,
      $h(z)=h(\pi(z),\phi_i(\pi(z)))
      =
      0$. 
      Note that we didn't say $z=(\pi(z),\phi_i(\pi(z))$, 
      but the value of $h$ is determined by the first $d+1$ coordinates.
      If $\pi(z)\in S$,
      since $W\setminus\pi^{-1}(S)$ is dense in $W$ and 
      zeroes of $p_{d+1}$ are zeroes of 
      $h$,
      by fact \ref{numofzeroes},
      every zero of $p_{d+1}$ in $W\cap\pi^{-1}(S)$ is also a zero of
      $h$, of at least the same multiplicity. 
      (i.e. If $z$ is a zero of $p_{d+1}$ of multiplicity $m$, 
      then $z$ is a zero of $h$ of multiplicity $\geq m$.)
      \\
      \indent
      Also by fact \ref{numofzeroes},
      for all $z'\in \pi(W)$,
      \begin{equation*}
\text{
     $|\{z\in W: \pi(z)=z' $ and $p_{d
     +1}(z)=0\}|\leq m$}.
   \end{equation*}
   Define the function
   $u$ on $W\setminus Z(p_{d+1})$,
   where $Z(p_{d+1})$ means the zero set of $p_{d+1}$, by
   \begin{equation*}
       u(z)=\dfrac{h(z)}{p_{d+1}(z)}.
   \end{equation*}
   Define a function $\overline{u}$ on $W$ by
   \begin{equation*}
       \overline{u}(z)=u(z) \text{ if
       $z\notin Z(p_{d+1})$};
   \end{equation*}
   \begin{equation*}
       \overline{u}(z)=\lim_{w\rightarrow z} u(w) \text{ if $z\in Z(p_{d+1})$}.
   \end{equation*}
   By \cite[Lemma~2.42]{peterzil2001expansions},
   for any $z'\in\pi(W)$,
   the function 
   $\overline{u}_{z'}(y):=
   \lim_{y'\rightarrow y}\dfrac{h(z',y)}{p_{d+1}(z',y)}$
   is  well-defined and $\mathcal{K}$-holomorphic, 
   since it has isolated singularities only.
   By \cite[Theorem~2.7]{peterzil2003expansions},
   $\overline{u}$ is continuous on $W$.
   By \cite[Theorem~2.14]{peterzil2003expansions},
   since $dim_{\mathcal{R}} Z(p_{d+1})\leq 2n-1$,
   $\overline{u}$ is $\mathcal{K}$-holomorphic on $W$.
   Since $h=\overline{u}\cdot p_{d+1}$
 on $W\setminus\pi^{-1}(S)$, 
 which is dense in $W$,
 $h=\overline{u}\cdot p_{d+1}$ on $W$.
       \end{proof}
      \begin{claim}(type version of \cite[Claim~11.8.]{peterzil2008complex})
          For every $g(z)\in\mathcal{O}_p$, 
          there is $h(z',u_{d+1},...,u_n)\in\mathcal{O}_{d,p}[\Bar{u}]$ of degree 
          less than $m$ in each variable $u_{d+1},...,u_n$ such that 
          $g(z)$ is equivalent to $h(z',z_{d+1},...,z_n)$ modulo $J_p$. 
          \end{claim}
          \begin{proof}
          Given $g(z)\in\mathcal{O}_p$ defined on some definable open 
          $U\subseteq\mathcal{K}^n$ with $p\in U$,
          may assume $\pi(U)\subseteq\pi(A)$ as in Claim \ref{claimdivisible}.
          Since $\forall z'\in\pi(U)$,
$p_{n}(z',w)$ has $m$ zeroes,
          by theorem \ref{tpvwdiv}, 
          there exist $k\in\mathbb{N}$ with $k\leq m$, a definable open 
          set $V\subseteq U$ with $p\in V$ and unique 
\begin{equation*}
   q(z',u_{d+1},...,u_n)\in\mathcal{O}_{V,n},
\end{equation*}
\begin{equation*}
R_0(z',u_{d+1},...,u_{n-1}),...,R_{k-1}(z',u_{d+1},...,u_{n-1})\in\mathcal{O}_{V,n-1}
\end{equation*}
such that 
\begin{equation*}
    g(z',u_{d+1},...,u_n)=q(z',u_{d+1},...,u_n)p_{n}(z',u_n)+
    \end{equation*}
    \begin{equation*}
      R_{k-1}(z',u_{d+1},...,u_{n-1})u_n^{k-1}+...+R_1(z',u_{d+1},...,u_{n-1})u_n+R_0(z',u_{d+1},...,u_{n-1}) 
    \end{equation*}
    on $V$.
Apply theorem \ref{tpvwdiv}. to $R_0,...,R_{k-1}$ by dividing $p_{n-1}(z',u_{n-1})$. 
Repeat this as in \cite[Claim 11.8.]{peterzil2008complex} and we get the conclusion.
          \end{proof}
     \begin{claim}\label{tpv11.9}
     (type version of \cite[Claim~11.9.]{peterzil2008complex})
         For every $g(z)\in\mathcal{O}_p$,  there is $q(z',u)\in\mathcal{O}_{d,p}[u]$ such that $\delta^sg(z)$ is equivalent modulo $J_p$ to $q(z',z_{d+1})$.
     \end{claim}
     \begin{proof}
      The same as in \cite[Claim 11.9.]{peterzil2008complex}.   
     \end{proof}
     \indent 
     We get Theorem \ref{tpvimmstp}. as in \cite[Theorem 11.2.]{peterzil2008complex} using Claim \ref{tpv11.9} and Claim \ref{claimdivisible}.
          
     \end{proof}
 \begin{theorem} (type version of \cite[Theorem~11.3.]{peterzil2008complex})
 \label{typever11.3}
 Assume that $A$ is a $\mathcal{K}$-analytic subset of $U\subseteq\mathcal{K}^n$ and assume that $G_1,...,G_t$ are $\mathcal{K}$-holomorphic maps from $A$ into $\mathcal{K}^N$.
 Then we can write $A$ as a union of finitely many relatively open sets $A_1,...,A_m$ such that on each $A_i$ the following holds:\\
 \indent
 There are finitely many tuples of $\mathcal{K}$-holomorphic functions on $A_i$, 
 \begin{equation*}
\{(H_{j,1},...,H_{j,t}):j=1,...,k\}, k=k(i),
 \end{equation*}
 with the property that for every $p\in S_n(\mathcal{K})$ with $p\in A_i$, 
 the module $R_p(g_1,...,g_t)$ equals its submodule generated by $\{(h_{j,1},...,h_{j,t}):j=1,...,k\}$ 
 over $\mathcal{O}_p$ (where $g_i$ and $h_{i,j}$ are the germs of $G_i$ and $H_{i,j}$ at $p$, resp).
 
 \end{theorem}
 \begin{proof}
 Induction on $n$. When $n=0$, $(G_1,..,G_t)$ can be considered as a vector in $\mathcal{K}^t$ and $\{(\phi_1,...,\phi_t):U\rightarrow\mathcal{K}^t:\phi_1G_1+...+\phi_tG_t=0\}$ is a vector subspace of $\mathcal{K}^t$.\\
 \indent 
 Assume true for $n-1$ and prove for $n$. Consider first $N=1$. 
 As in \cite[Claim~1 in Theorem~11.3.]{peterzil2008complex}, 
 may assume $G_1,...,G_t$ are Weierstrass polynomials $\omega_1(z',z_n),....,\omega_t(z',z_n):V\times\mathcal{K}=\pi(U)\times\mathcal{K}\rightarrow\mathcal{K}$,
 where $\pi$ is the projection onto the first $n-1$ coordinates.\\
 \indent
 We reduce to the case where $(\phi_1,...,\phi_t)$ are polynomials as in \cite[Claim~2 in Theorem~11.3.]{peterzil2008complex}: \\
 \begin{claim}(type version of \cite[Claim~2 in Theorem~11.3.]{peterzil2008complex})
 For $p\in S_n(\mathcal{K})$ with $p\in U'$, $U'\subseteq U$ definable and open, let $\phi=(\phi_1,...,\phi_s)\in R_{U'}(\omega)$. 
 \\
 \indent
 Then there are tuples $\{\psi^i=(\psi^i_1,...,\psi^i_s):i=1,...,t\}$, where each $\psi^i_j$ is a polynomial in $z_n$ of degree $\leq m$ over $\mathcal{O}_{V''}$ ($V''=\pi(U'')$, where $U''\subseteq U'$ is definable open and $p\in U''$), such that $\phi$ is in the module generated by $\psi^1,...,\psi^t$ over $\mathcal{O}_{U''}$.       
 \end{claim}
 \begin{proof}
 Let $V=\pi(U)$. 
 By lemma \ref{tubu}, 
 there exist $k\in\mathbb{N}$ and
 $U''\subseteq U'$ a definable open neighborhood of $p$ 
        such that for any $z\in \pi(U'')$,
        there are exactly $k$ zeroes of
        $\omega_1(z,-)$ in $U''_y$ (counting multiplicity).
     Let $h_1,...,h_k$
        be definable functions
        such that
        for all $z\in \pi(U'')$,
        $\{h_1(z),...,h_k(z)\}$ list all zeroes
        of $\omega_1(z,-)$ in $U_y''$.
 Take 
 \begin{equation*}
 \text{
     $\omega_1'(z',y)=\underset{i=1}{\overset{k}{\prod}} (y-h_i(z'))$.}
      \end{equation*}
$\omega_1'$ is holomorphic
by the same argument as in \cite[Theorem~2.20]{peterzil2003expansions}.
 By theorem \ref{tpvwdiv}, there exist definable unique $q_i(z',z_n)\in\mathcal{O}_{U''}$, polynomials
 $r_i(z',z_n)\in\mathcal{O}_{U''}$ of degree $<k$ such that 
 $\phi_i=q_i\omega_1'+r_i$ on $U''$. Since $\forall (z',y)\in U''$, 
 $\omega_1'(z',y)=0$ implies $\omega_1(z',y)=0$ and for any such zero 
 $(z',y)\in U''$, the degree of the zero $y$ in $\omega_1'(z'-)$ is the same as 
 the degree of that zero in $\omega_1(z'-)$, by the same construction as in Claim \ref{claimdivisible},
 there is $u'\in\mathcal{O}_{U''}$ such that $\omega_1=u'\omega_1'$. 
 Moreover, $\forall (z',y)\in U''$, $u'\neq 0$. 
 Define
 \begin{equation*}
     \psi^2=(-\omega_2,\omega_1,0,...,0),
         ..., 
      \end{equation*}
     \begin{equation*}
     \psi^s=(-\omega_s,0,...,0,\omega_1), 
   \end{equation*}
   \begin{equation*}
       \psi=\phi_1+q_2\omega_2+...+q_s\omega_s,
       \end{equation*}
       \begin{equation*}
       P=-(r_2\omega_2+...+r_s\omega_s), 
        \end{equation*}
       \begin{equation*}
       \text{$Q=$ 
 a $\mathcal{K}$-holomorphic polynomial in variable $z_n$ over $\pi(U'')$ extending $u'\psi$, }
  \end{equation*}
 \begin{equation*}
 \text{
 $Q_i=$ a $\mathcal{K}$-holomorphic polynomial in variable $z_n$ extending 
 $u'r_i$ for $i\in\{2,...,s\}$, }
 \end{equation*}
 and 
 \begin{equation*}
     \text{$\psi_1=(Q,Q_2,...,Q_s)$}.
      \end{equation*} 
      As in 
 \cite[Chapter~8, Theorem~8B]{whitney1972complex},
 \begin{equation*}
\text{ $\phi=\underset{i=2}{\overset{s}{\sum}}q_i\psi_i+(1/u')(Q,Q_2...,Q_s)$ in $U''$.} 
 \end{equation*}
 (For the same reason as in \cite[Chapter~8, Theorem~8B]{whitney1972complex}, it's important to 
 use $\omega_1'$ instead of $\omega_1$ because there might be zeros of $\omega_1$ outside the domain of $\phi$, 
 i.e. outside $U'$. When there exist zeros of $\omega_1$ outside the domain of $\phi$,
 we cannot use the equation $\psi\omega_1=P$ to conclude that zeros 
 of $\omega_1$ are also zeros of $P$, and hence cannot extend $\psi$ as $P/\omega_1$.).    
 \end{proof}
 The rest is the same argument as in \cite[Theorem~11.3.]{peterzil2008complex}
 \end{proof}
 \begin{theorem} \label{typecoh} (type version of \cite[Theorem~11.1.]{peterzil2008complex})
 Let $M$ be a $\mathcal{K}$-manifold and $A\subseteq M$ a $\mathcal{K}$-analytic subset of $M$. 
 Then there are finitely many open sets $V_1,...,V_k$ whose union covers $M$ and for each $i=1,...,k$ there are finitely many $\mathcal{K}$-holomorphic functions $f_{i,1},...,f_{i,m_i}$ in 
 $\mathcal{I}_{V_i}(A)$, such that for every $p\in S_n(\mathcal{K})$ with $p\in V_i$ the functions 
 $f_{i,1},...,f_{i,m_i}$ generate the ideal $\mathcal{I}(A)_p$ in $\mathcal{O}_p$.\\
 \indent
 Moreover, the $V_i$'s and the $f_{i,j}$'s are all definable over the same parameters defining $M$ and $A$.
     
 \end{theorem}
 \begin{proof}
 The same as in \cite{peterzil2008complex} using Theorem \ref{tpvimmstp} and Theorem \ref{typever11.3}.    
 \end{proof} 
 \begin{theorem}\label{mainthm} (another proof of \cite[Theorem 2.21]{bakker2022minimal})
     The definable structure sheaf $\mathcal{O}_{\mathcal{K}^n}$ of $\mathcal{K}^n$ is a coherent
     $\mathcal{O}_{\mathcal{K}^n}$-module as a sheaf on the site $\underline{\mathcal{K}^n}$.
 \end{theorem}
 \begin{proof}
     $\mathcal{O}_{\mathcal{K}^n}$ is a generated by $1$ as an $\mathcal{O}_{\mathcal{K}^n}$-module. \\
     \indent
     It suffices to show that given any definable open $U \subseteq \mathcal{K}^n$ and any $\mathcal{O}_U$-module homomorphism $\varphi: \mathcal{O}^m_U \rightarrow \mathcal{O}_U$, the kernel of $\varphi$ is of finite type. i.e. we want a finite definable cover $U_i$ of $U$ and surjections $\mathcal{O}^n_{U_i} \twoheadrightarrow $ $(ker\varphi)_{_{U_i}}$ for some positive
integer $n$ on each of those open sets. Let $G_1,...,G_m$ be definable $\mathcal{K}$-holomophic functions from $U$ to $\mathcal{K}$ such that $\varphi(e_j)=G_j$ for all $e_j$ in the canonical basis $\{e_1,...,e_m\}$ of $\mathcal{O}^m_U$. By theorem \ref{typever11.3}, $U$ is a union of finitely many definable open sets $U_1,...,U_l$ such that on each $U_i$ the following holds:\\
 \indent
  There are finitely many tuples of $\mathcal{K}$-holomorphic functions on $U_i$, $\{(H_{j,1},...,H_{j,m}):j=1,...,k\}$, $k=k(i)$, with the property that for every $p\in S_n(\mathcal{K})$ with $p\in U_i$, the module $R_p(g_1,...,g_m)$ equals its submodule generated by $\{(h_{j,1},...,h_{j,m}):j=1,...,k\}$ over $\mathcal{O}_p$ (where $g_i$ and $h_{i,j}$ are the germs of $G_i$ and $H_{i,j}$ at $p$, resp).\\
  \indent
  Hence, given $U_i$ as above and a definable open $V\subseteq U_i$, fix any $s\in ker(\varphi)_{U_i}(V)$. For each $p\in S_n(\mathcal{K})$ with $p\in V$, since $s_p\in ker(\varphi)_p=R_p(g_1,...,g_m)$, there is a definable open $V_p$ with $p\in V_p\subseteq V$ such that $s|_{V_p}$ is generated by $\{(H_{j,1}|_{V_p},...,H_{j,m}|_{V_p}):j=1,...,k\}$, $k=k(i)$ on $V_p$.  By compactness of $S_n(\mathcal{K})$, there exist finitely many $p_1,...,p_s$ such that $V=V_{p_1}\cup...\cup V_{p_s}$ and on each of these $V_{p_\alpha}$'s, $s|_{V_{p_\alpha}}$ is generated by $\{(H_{j,1}|_{V_{p_{\alpha}}},...,H_{j,t}|_{V_{p_{\alpha}}}):j=1,...,k\}$, $k=k(i)$. By definition \ref{injsur}, the morphism $\mathcal{O}_{U_i}^{k(i)}\rightarrow ker(\varphi)_{U_i}$ given by mapping the canonical basis to $\{(H_{j,1}|_{U_{i}},...,H_{j,t}|_{U_{i}}):j=1,...,k\}$ is surjective. Hence $ker(\varphi)$ is of finite type and $\mathcal{O}_{\mathcal{K}^n}$ is a coherent $\mathcal{O}_{\mathcal{K}^n}$-module.
  
  \end{proof}
  \section{Remark}\label{rem}
  \indent
    
    Let $X\subseteq\mathcal{K}^n$ be definable open. 
    Let $Sh^{\mathcal{O}_X}(\Tilde{X})$ be the category of usual sheaves $\mathcal{O}_X$-modules on $\Tilde{X}$.
    Let $Sh_{dtop}^{\mathcal{O}_X}(X)$ be the category of sheaves of $\mathcal{O}_X$-modules on $\underline{X}$ as an o-minimal site.
  We show that $Sh^{\mathcal{O}_X}(\Tilde{X})$ and $Sh_{dtop}^{\mathcal{O}_X}(X)$ are isomorphic categories, 
  and 
   the surjective maps are exactly the epimorphisms in both categories.  
   Hence,  
from a category-theoretic perspective, theorem \ref{typecoh} immediately implies
theorem \ref{mainthm}.
  \subsection{Sheafification}
  \begin{definition}
  \cite[Part~1, Chapter~7, Section~10]{stacks-project} Let $X\subseteq \mathcal{K}^n$ be a definable open set. Let $\mathcal{F}$ be a presheaf of $\mathcal{O}_X$-modules, and let $U=\{U_i\}_{i=1}^{k}$ be a covering of $U$. Let us use the notation $\mathcal{F}(U)$ to indicate the equalizer 
  \begin{equation*}
      H^0(U, \mathcal{F}) = \{(s_i)_{i\in \{1,...,k\}} \in\prod_i \mathcal{F}(U_i) : s_i|_{U_i\cap U_j} = s_j|_{U_i\cap U_j} \forall i, j \in \{1,...,k\}\}.
      \end{equation*}
      \end{definition}
      \begin{definition}\label{sheafidef}
      \cite[Part~1, Chapter~7, Section~10]{stacks-project}
      Let $X\subseteq\mathcal{K}^n$ be a definable open set. 
      For $U\subseteq X$ be definable open, 
      let $J_U$ be the set of all finite definable open coverings of $U$. 
      Define $\leq$ on $J_U$ by $U\leq V$ if $V$ is a refinement of $U$. $(J_{U},\leq)$ is a directed set. 
      For $V=\{V_j\}_{j=1}^l$ a refinement of $U=\{U_i\}_{i=1}^k$, 
      fix a function $\alpha:[l]\rightarrow[k]$ such that $V_j\subseteq U_{\alpha(j)}$. 
      Define $\mu_{U,V}:H^0(U,\mathcal{F})\rightarrow H^0(V,\mathcal{F})$ by 
      $\mu_{U,V}((s_i:i=1,...k))=(s_{\alpha(j)}|_{V_j}:j=1,...,l)$. 
      (In fact, $\mu_{U,V}$ is independent of the choice of 
      $\alpha$: if $V_j\subseteq U_i$ and $V_j\subseteq U_{i'}$, 
      then $V_j\subseteq U_i\cap U_{i'}$ 
      and $s_i|_{V_j}=s_i|_{U_i\cap U_{i'}}|_{V_j}=s_{i'}|_{U_i\cap 
      U_{i'}}|_{V_j}=s_{i'}|_{V_j}$.) ($H^0(U,\mathcal{F}),\mu_{U,V}$) is a directed system. 
      Define the presheaf $\mathcal{F}^+$ by
         \begin{equation*}
              \mathcal{F}^+(U) = \coprod_{U\in J_U}H^0(U,\mathcal{F})/ \sim
\end{equation*}
\text{where for} $s \in H^0(U,\mathcal{F})$ and $s' \in H^0(V,\mathcal{F})$ \text{we have} $s \sim s' \iff \mu_{U,\mathcal{W}}(s) = \mu_{V,\mathcal{W}}(s')$ \text{for some} $\mathcal{W} \geq U, V$.\\
\indent
For a presheaf $\mathcal{F}$, define the canonical map $\tau:\mathcal{F}\rightarrow\mathcal{F}^{+}$ by $\mathcal{F}(U)\rightarrow\mathcal{F}^{+}(U):s\mapsto (s)/\sim$. 
      \end{definition}
      \begin{fact}\label{sheafithm}
      \cite[Part~1, Chapter~7, Section~10, Theorem~10.10]{stacks-project}
      \begin{enumerate}
          \item [(1)]
      The presheaf $\mathcal{F}^{+}$ is separated.
        \item [(2)]
      If $\mathcal{F}$ is separated, then $\mathcal{F}^{+}$ is a sheaf and the map of presheaves $\mathcal{F}\rightarrow\mathcal{F}^{+}$ is injective.
        \item [(3)] 
        If $\mathcal{F}$ is a sheaf, then $\mathcal{F}\rightarrow\mathcal{F}^{+}$ is an isomorphism.
       \item [(4)]
      The presheaf $\mathcal{F}^{++}$ is always a sheaf.
          \end{enumerate}
      \end{fact}  
      \begin{definition}\cite[Part~1, Chapter~7, Section~10, Definition~10.11]{stacks-project}
          The sheaf $\mathcal{F}^{\#} := \mathcal{F}^{++}$ together with the canonical map $\tau^{\#}=\tau^{+}\circ\tau:\mathcal{F} \rightarrow \mathcal{F}^{+}\rightarrow\mathcal{F}^{\#}$ is called the \emph{sheaf associated to $\mathcal{F}$}.
      \end{definition}
      \subsection{Epimorphism}
     \indent
     
      Following proofs in \cite{stacks-project},
      we show that surjective maps and epimorphisms coincide in the categories $Sh_{dtop}^{\mathcal{O}_X}(X)$ and $Sh^{\mathcal{O}_X}(\Tilde{X})$
      \begin{lemma}\label{surepi}
          \cite[Part~1, Chapter~7, Lemma~11.2.]{stacks-project} The surjective maps defined in definition \ref{injsur} are exactly the epimorphisms of the category $Sh_{dtop}^{\mathcal{O}_X}(X)$. 
      \end{lemma}
    \begin{proof}
        Let $\varphi:\mathcal{F}\rightarrow\mathcal{G}$ be an epimorphism between $\mathcal{F}$, $\mathcal{G}$ which are $\mathcal{O}_X$-modules on the o-minimal site $\underline{X}$. Consider the presheaf $\mathcal{H}$ defined by $\mathcal{H}(U)=\mathcal{G}(U)\oplus\mathcal{G}(U)/S(U)$ where $S(U)$ is the $\mathcal{O}(U)$-submodule $\{(y,z)\in\mathcal{G}(U)\oplus\mathcal{G}(U):\exists x\in \mathcal{F}(U)$ $\varphi_{U}(x)=y,z=-y\}$ for $U\subseteq X$ definable open (i.e. the pushout). \\
        \indent
        As in \cite[Part~1, Chapter~7, Section~3, Lemma 3.2.]{stacks-project}, consider the presheaf morphisms $i_1,i_2:\mathcal{G}\rightarrow\mathcal{H}$ defined by $i_{1_U}(x)=(x,0)/S(U)$ and $i_{2_U}(x)=(0,x)/S(U)$. Let $i_1'=\tau^{+}\circ\tau\circ i_1:\mathcal{G}\rightarrow\mathcal{H}^{\#}$, $i_2'=\tau^{+}\circ\tau\circ i_2:\mathcal{G}\rightarrow\mathcal{H}^{\#}$. Then $i_1'$, $i_2'$ are morphisms in $Sh_{dtop}^{\mathcal{O}(X)}(X)$. Since $i_1\circ\varphi=i_2\circ\varphi$ as presheaf morphisms by definition, $i_1'\circ\varphi=i_2'\circ\varphi$ as sheaf morphisms. Since $\varphi$ is an epimorphism, $i_1'=i_2'$. Fix $U\subseteq X$ definable open and $y\in \mathcal{G}(U)$. Since $i_1'(y)=i_2'(y)$ and $\tau^{+}$ is injective, by fact \ref{sheafithm} (1), (2), $\tau(i_1(y))=\tau(i_2(y))$. By definition \ref{sheafidef}, there exists a finite definable open covering $\{U_i\}_{i=1}^k$ of $U$ such that $(y|_{U_i},0|_{U_i})/S(U_i)=(0|_{U_i},y|_{U_i})/S(U_i)$ for all $i=1,...,k$. By definition of $\mathcal{H}$, for each $i\in\{1,...,k\}$, there exists $x_i\in\mathcal{F}(U_i)$ such that $\varphi_{U_i}(x_i)=y|_{U_i}$. Hence $\varphi$ is a surjective morphism.\\
        \indent
        The other direction is just checking definitions.

    \end{proof}  
   We have a similar result for 
     $Sh^{\mathcal{O}_X}(\Tilde{X})$,
     the category of (classical) sheaves on $\Tilde{X}$ as a topological space. 
    \begin{lemma}
        The surjective maps (i.e. surjective at the stalks) are exactly the epimorphisms of the category $Sh^{\mathcal{O}_X}(\Tilde{X})$. 
    \end{lemma}
    \begin{proof}
        The same as in the proof lemma \ref{surepi} using the usual sheafification of sheaves.
    \end{proof}
    \begin{proposition}
        $Sh^{\mathcal{O}_X}(\Tilde{X})$ and $Sh_{dtop}^{\mathcal{O}_X}(X)$ are isomorphic categories. 
    \end{proposition}
    \begin{proof}
        The same proof as in \cite[Proposition~3.2]{edmundo2006sheaf}
    \end{proof}
    
    Another proof of theorem \ref{mainthm}: Let $\iota:Sh^{\mathcal{O}_X}(\Tilde{X})\rightarrow Sh_{dtop}^{\mathcal{O}_X}(X)$ be an isomorphism. 
    Let  $U \subseteq \mathcal{K}^n$ be definable open and $\varphi: \mathcal{O}^m_U \rightarrow \mathcal{O}_U$ a $\mathcal{O}_U$-module homomorphism.
    By theorem \ref{typever11.3}, 
    there exists a finite definable open covering $\{U_i\}_{i=1}^k$ of $U$ such that for some $l\in\mathbb{N}$ and for each $i\in\{1,...,k\}$,
    there exists $\psi_i:\mathcal{O}_{\Tilde{U}_i}^l\twoheadrightarrow ker(\iota^{-1}(\varphi))_{\Tilde{U}_i}$.
    Since surjective morphisms are epimorphisms in $Sh^{\mathcal{O}_X}(\Tilde{X})$,
    $\iota(\psi_i):\mathcal{O}_{U_i}^l\rightarrow ker(\varphi)_{U_i}$ is an epimorphism and hence a surjective morphism by lemma \ref{surepi}.
  \bibliographystyle{alpha}
\bibliography{ref}

\end{document}